\documentclass[11pt]{amsart}
\usepackage{amssymb,latexsym}

\newdimen\AAdi%
\newbox\AAbo%
\def\AAk#1#2{\s_etbox\AAbo=\hbox{#2}\AAdi=\wd\AAbo\kern#1\AAdi{}}%
\def\AAr#1#2#3{\s_etbox\AAbo=\hbox{#2}\AAdi=\ht\AAbo\raise#1\AAdi\hbox{#3}}%
\font\tenmsb=msbm10 at 12pt
\font\sevenmsb=msbm7 at 8pt
\font\fivemsb=msbm5 at 6pt
\newfam\msbfam
\textfont\msbfam=\tenmsb
\scriptfont\msbfam=\sevenmsb
\scriptscriptfont\msbfam=\fivemsb
\def\Bbb#1{{\tenmsb\fam\msbfam#1}}

\begin{document}

\newtheorem{thm}{Theorem}
\newtheorem{lem}{Lemma}
\newtheorem{cor}{Corollary}
\newtheorem{rem}{Remark}
\newtheorem{pro}{Proposition}
\newtheorem{defi}{Definition}
\newcommand{\noi}{\noindent}
\newcommand{\dis}{\displaystyle}
\newcommand{\mint}{-\!\!\!\!\!\!\int}
\newcommand{\ba}{\begin{array}}
\newcommand{\ea}{\end{array}}

\def \bx{\hspace{2.5mm}\rule{2.5mm}{2.5mm}} \def \vs{\vspace*{0.2cm}}
\def\hs{\hspace*{0.6cm}}
\def \ds{\displaystyle}
\def \p{\partial}
\def \O{\Omega}
\def \o{\omega}
\def \b{\beta}
\def \m{\mu}
\def \l{\lambda}
\def\L{\Lambda}
\def \ul{u_\lambda}
\def \D{\Delta}
\def \d{\delta}
\def \s{\sigma}
\def \e{\varepsilon}
\def \a{\alpha}
\def \tf{\tilde{f}}
\def\cqfd{%
\mbox{ }%
\nolinebreak%
\hfill%
\rule{2mm} {2mm}%
\medbreak%
\par%
}
\def \pr {\noindent {\it Proof.} }
\def \rmk {\noindent {\it Remark} }
\def \esp {\hspace{4mm}}
\def \dsp {\hspace{2mm}}
\def \ssp {\hspace{1mm}}

\def \u{u_+^{p^*}}
\def \ui{(u_+)^{p^*+1}}
\def \ul{(u^k)_+^{p^*}}
\def \energy{\int_{\R^n}\u }
\def \sk{\s_k}
\def \mo{\mu_k}
\def\cal{\mathcal}
\def \I{{\cal I}}
\def \J{{\cal J}}
\def \K{{\cal K}}
\def \OM{\overline{M}}

\def\fk{{{\cal F}}_k}
\def\M1{{{\cal M}}_1}
\def\Fk{{\cal F}_k}
\def\Fl{{\cal F}_l}
\def\FF{\cal F}
\def\Gk{{\Gamma_k^+}}
\def\n{\nabla}
\def\uuu{{\n ^2 u+du\otimes du-\frac {|\n u|^2} 2 g_0+S_{g_0}}}
\def\uuug{{\n ^2 u+du\otimes du-\frac {|\n u|^2} 2 g+S_{g}}}
\def\sku{\sk\left(\uuu\right)}
\def\qed{\cqfd}
\def\vvv{{\frac{\n ^2 v} v -\frac {|\n v|^2} {2v^2} g_0+S_{g_0}}}
\def\vvs{{\frac{\n ^2 \tilde v} {\tilde v}
 -\frac {|\n \tilde v|^2} {2\tilde v^2} g_{S^n}+S_{g_{S^n}}}}
\def\skv{\sk\left(\vvv\right)}
\def\tr{\hbox{tr}}
\def\pO{\partial \Omega}
\def\dist{\hbox{dist}}
\def\RR{\mathcal R}
\def\R{\Bbb R}
\def\C{\Bbb C}
\def\B{\Bbb B}
\def\N{\Bbb N}
\def\Q{\Bbb Q}
\def\Z{\Bbb Z}
\def\PP{\Bbb P}
\def\EE{\Bbb E}
\def\F{\Bbb F}
\def\G{\Bbb G}
\def\H{\Bbb H}
\def\SS{\mathcal S}
\def\S{\Bbb S}

\def\be{\begin{equation}}
\def\ee{\end{equation}}
\def\ba{\begin{array}}
\def\ea{\end{array}}

\def\lcf{{locally conformally flat} }

\def\Ric{{\rm Ric}}
\def\div{{\rm div}}

\date{}

\title[Conformal deformations of  Ricci tensor]
{Conformal deformations of the smallest eigenvalue of the Ricci
tensor}
\author{Pengfei Guan}
\address{Department of Mathematics\\
McGill University\\
Montreal, H3A 2K6, Canada.} \email{guan@math.mcgill.ca}
\thanks{Research of the first author was supported in part by
an NSERC Discovery Grant.}
\author{Guofang Wang}
\address{Max Planck Institute for Mathematics in
the Sciences\\ Inselstr. 22-26, 04103 Leipzig, Germany}
\email{gwang@mis.mpg.de}

\begin{abstract}
We consider deformations of metrics in a given conformal class
such that the smallest eigenvalue of the Ricci tensor to be a
constant. It is related to the notion of minimal volumes in
comparison geometry. Such a metric with the smallest eigenvalue of
the Ricci tensor to be a constant is an extremal metric of volume
in a suitable sense in the conformal class. The problem is reduced
to solve a Pucci type equation with respect to the Schouten
tensor. We establish a local gradient estimate for this type of
conformally invariant fully nonlinear uniform elliptic equations.
Combining it with the theory of fully nonlinear equations, we
establish the existence of solutions for this equation.

\end{abstract}
\maketitle
\section{introduction}

\medskip

The Ricci curvature tensor of a Riemannian metric plays an
important role in comparison geometry for Riemannian manifolds, in
particular the lower bounds of Ricci curvature. See for instance
\cite{Bishop, Gro2, petersen2}. In this paper, we
are interested in conformal deformations of the smallest eigenvalue
of the Ricci tensor. Let $(M^n,g_0)$ be an $n$-dimensional compact
Riemannian manifold and $[g_0]$  its conformal class. And let
$Ric_g$ and $R_g$  be the Ricci curvature tenser and the scalar
tensor of a metric $g$ respectively. Define $\min Ric_g(x)$ the
smallest eigenvalue of $g^{-1}\cdot Ric_g$ at $x\in M$. Our
problem is to find a conformal metric $g=e^{-2u}g_0$ such that
\begin{equation}
\label{eq1}
\min Ric_g(x)=constant.
\end{equation}

It turns out that the problem is equivalent to solving an
interesting fully nonlinear uniform elliptic equation. First we
recall that the Schouten tensor of the metric $g$ is defined by
\[S_g=\frac{1}{n-2} \left( Ric_g - \frac{R_g}{2(n-1)}g\right).\]
Let $\L=(\l_1, \l_2,\cdots,\l_n)\in \R^n$. Assume that $\l_1\le
\l_2\le\cdots \le \l_n$. For an integer $1\le p\le n-1$, define a
function $G_p:\R^n\to \R$ by
\[G_p(\L)=
(n-p)\sum _{i\le p}\l_{i}+p\sum _{i>p}\l_{i}.
\]
For a symmetric matrix $A$, $G_p(A)=G_p(\L)$, where $\L$ is the set of
eigenvalues of $A$. It is easy to check $\min Ric=
G_1(g^{-1}\cdot S_g)$.

We may also ask if there is a conformal metric with a constant
$W_p(g):=G_p(g^{-1}\cdot S_g)$. $W_p(g)$ is also an interesting
geometric object, which will be called {\it $p$-Weitzenb\"ock
curvature}, for it arises from the Weitzenb\"ock formula for
$p$-forms in a \lcf manifold. See the Appendix.  We will
consider the following general equation
\begin{equation}\label{eq4.0}W_p(g)(x)=f(x),
\quad\forall x\in M.\end{equation}

We first treat the case when the background metric has negative
curvature, i.e., $W_p(g_0)<0$.
\begin{thm}\label{thm2} Let $(M,g_0)$ be a compact Riemannian manifold and
$1\le p\le n$. Suppose that $W_p( g_0)(x) < 0$ for any $ x\in M$,
then there is a unique $C^{2,\a}$ metric $g\in [g_0]$ for some
$\alpha>0$ such that $W_p(g)(x)=-1, \quad\forall x\in M.$
\end{thm}

A geometric consequence of Theorem \ref{thm2} is the existence of
an extremal metric in the conformal class with minimal volume.
Although (\ref{eq4.0}) has no variational structure in general, a
solution of (\ref{eq4.0}) does achieve the minimum of the minimal
volume in a conformal class in this case. The following is a
simple consequence, see section 4 for other related results.

\begin{cor}\label{ivcomp}
Suppose that $\min Ric_{g_0}(x)< 0$ for any $x\in M$. Then there
is a unique conformal metric $g^*\in [g_0]$ such that
$vol({g^*})=\min {vol(g)}$, where minimum is taken over all $g\in
[g_0]$ with $\min Ric_{g}(x)\ge -1$. The extremal metric $g^*$ is
characterized by a unique solution to equation $\min
Ric_{g^*}(x)=-1, \forall x\in M$.
\end{cor}

This Corollary is related to the minimal volumes considered by Gromov \cite{G}
and LeBrun \cite{Le}
\medskip

We now turn to the case of the positive Ricci curvature.

\begin{thm}\label{thm3} Let $(M,g_0)$ be a compact Riemannian manifold
with $Ric_{g_0}>0$.  Then there is a conformal metric $ g\in
[g_0], g\in C^{2,\alpha}(M)$ for some $\alpha>0$ such that
$\min Ric_{g}(x)=n-1$ for all $x\in M$. \end{thm}

The positivity of $p$-Weitzenb\"ock curvature for $1<p\le n/2$ plays
an important role in the investigation of the topological structure of
locally conformally flat manifolds in \cite{GLW} and the Appendix.

\begin{thm}\label{thm5}
Let $(M,g_0)$ be an $n$-dimensional smooth compact locally
conformally flat manifold with $W_p(g_0)>0$ and $p\le n/2$. If
$(M,g_0)$ is not conformally equivalent to the standard
$n$-sphere, then there exists $ g\in [g_0], g\in C^{2,\alpha}(M)$
for some $\alpha>0$ such that
\begin{equation}\label{eqb1}W_p(g)(x)=1, \quad\forall x\in M.\end{equation}
Furthermore, the solution space is
compact. That is, there is positive constant $C>0$, such that
\begin{equation}\label{glwgradest}
\|u\|_{C^{2,\alpha}} \leq C\end{equation} for any $C^{2,\alpha}$
solution $g=e^{-2u}g_0$ of (\ref{eqb1}).
\end{thm}

One of the interesting aspect of this paper is that our geometric
problem is deduced to a Pucci type fully nonlinear uniformly
elliptic equations with respect to the Schouten tensor. These
equations have been studied extensively in Euclidean domains, in
particular in connection to stochastic optimization. With the
breakthrough of the Krylov-Safonov's Harnack estimate \cite{KS}
for non-divergent elliptic equations, it followed the fundamental
Evans-Krylov theorem \cite{Evans, Krylov} on the H\"older
regularity of the second derivatives. The deep results of
Caffarelli \cite{CC} and Safonov \cite{Safonov1, Safonov} provide
further improvement of such H\"older estimates bypassing $C^2$
estimates. Caffarelli's perturbation theory is very important for
our application in conformal geometry. Especially some of our
geometric equations are in terms of convex operators, where $C^2$
estimates are difficult to obtain with standard approaches. The
core of Caffarelli's perturbation theory is the $C^2$ estimates
for simple model equations, which were established by him for the
concave and convex operators. To apply the Caffarelli perturbation
theory for our equations, we need one more crucial estimate, that
is a local gradient bound depending only on one side bound of
solutions. As in the case for nonlinear equations in conformal
geometry, this is a key step for the blow-up analysis. The Bishop-Gromov
volume comparison theorem also pays an important role in the proof of
Theorem \ref{thm5}.

\medskip

The paper is organized as follows. In Section 2, we discuss
briefly our uniformly elliptic fully nonlinear equations and
Caffarelli's results. We prove the local gradient estimates in
Section 3. In Section 4 we discuss the case of negative curvature
and the minimal volume in a conformal class. We prove Theorem
\ref{thm3} and Theorem \ref{thm5} in Section 5.

\medskip

{\it Acknowledgement.} We would like to thank Claude LeBrun for
 helpful conversations with the first author at AIM's workshop
on ``Conformal Structure in Geometry" in Palo Alto. When informed
our results on the conformal deformations of the smallest
eigenvalue of Ricci tensor, he pointed to us a possible
application to the minimal volumes. Corollary \ref{ivcomp} was
inspired by conversations with him.

\section{Fully nonlinear uniformly elliptic equations and Caffarelli's
estimates}

In this section, we will deduce our problem to fully nonlinear
uniformly elliptic elliptic equations.  We will make use of
Caffarelli's fundamental $W^{2,p}$ and $C^{2,\alpha}$ estimates in
\cite{CC} (see also Safonov's work on $C^{2,\alpha}$ estimates in
\cite{Safonov}). Some of these results have been subsequently
generalized to certain type of equations of form $F(\nabla^2 u,
\nabla u, x)=f(x)$ by L. Wang in \cite{LW} as well as for the case
of parabolic equations. As we will see that equation (\ref{eq4.0})
involves $\nabla u$ in a delicate way, we need certain appropriate
a priori estimates depending only on one side bound of $u$ (to be
more explicit, $e^{-2\inf u}$). This type of local gradient
estimates will be established in the next section.

Let $(M^n,g_0)$ be a compact, oriented Riemannian manifold of
dimension $n>2$. Let $[g_0]$ be the conformal class of $g_0$. For
any $g\in [g_0]$, we denote $Ric_g$, $S_g$ the Ricci tensor and
the Schouten tensor of the metric $g$ respectively. We write
\[\begin{array}{rcl}
\ds\vs {\mathcal
R}_g^s(x)&=&\ds\mbox{smallest eigenvalue of} \quad g^{-1}Ric_g(x),\\
\ds\vs{\mathcal R}_g^l(x)&=&\ds\mbox{largest eigenvalue of} \quad
g^{-1}Ric_g(x),\\
\ds\vs{\mathcal S}_g^s(x)&=&\ds\mbox{smallest
eigenvalue of} \quad g^{-1}S_g(x),\\
\ds {\mathcal
S}_g^l(x)&=&\ds\mbox{largest eigenvalue of} \quad g^{-1}S_g(x).
\end{array}\]
It is clear that $\RR^s_g$ and  $\RR^l_g$ do not depend on the
choice of coordinates. From the definition of the Schouten tensor
\[
S_g=\ds\vs \frac{1}{n-2}\left(Ric_g-\frac {R_g}{2(n-1)}\right),
\]
we have the following relations
\begin{equation}\label{eq2.1}\begin{array}{rcl}
\ds {\mathcal R}_g^s(x)&=& \ds\vs (n-2){\mathcal S}_g^s(x)+\tr(S_g(x)),\\
\ds {\mathcal R}_g^l(x)&=& \ds(n-2){\mathcal S}_g^l(x)+\tr(S_g(x)).\end{array}
\end{equation}

If $g=e^{-2u}g_0$, there is a transformation formula between two
Schouten tensors
\begin{equation}\label{transf} S_g=\nabla^2 u+\nabla u \otimes
\nabla u-\frac{|\nabla u|^2 g_0}{2}+S_{g_0},\end{equation} where
all covariant derivatives are with respect to $g_0$. (The same
convention will be used in the rest of this paper, unless it is
stated otherwise). Locally, denoting $S_{ij}$ the Schouten tensor
of $g_0$ under the frame, we have
\[(S_g)_{ij}=u_{ij}+u_iu_j-\frac{|\nabla u|^2}2
\delta_{ij}+S_{ij}.\]
\medskip

For any symmetric matrix $A$, we denote $\lambda_s(A)$ and
$\lambda_l(A)$ to be the smallest and largest eigenvalues of $A$
respectively. And we denote $\s_1(A)=\tr (A)$. We obtain the
equations for constant ${\mathcal R}_g^s$ and ${\mathcal R}_g^l$
respectively:
\begin{equation}\label{sR}
F_{*}(u)=: ((n-2)\lambda_s+\sigma_1)(\uuu)
=e^{-2u}{\mathcal R}_g^s,\end{equation}
\begin{equation}\label{lR}
F^{*}(u)=: ((n-2)\lambda_l+\sigma_1)(\uuu)=e^{-2u}{\mathcal R}_g^l.\end{equation}

Hence $F_{*}$ and $F^{*}$ are uniformly elliptic with ellipticity
constants  $1$ and $n-1$. It is clear that $F_{*}$ is concave and
$F^{*}$ is convex. We also note that $F^*$ and $F_*$ are
homogeneous of degree $1$.

\medskip

There are other similar fully nonlinear equations arising in the
Weitzenb\" ock formula for $p$-form on local conformally flat manifolds.
Let $\L=(\l_1, \l_2,\cdots,\l_n)\in \R^n$ be the set of eigenvalues of
a symmetric matrix $A$. For an integer $1\le p\le n/2$,
define a function
$G_p:\R^n\to \R$ by
\[G_p(A)=G_p(\L)=\min\{ (n-p)\sum _{k\le p}\l_{i_k}+p\sum _{k>p}\l_{i_k}\},\]
where $\min$ is over all permutations of $1,2,\cdots, n$. We
define $W_p(g)$ the  $p$-Weitzenb\"ock curvature of $g$ by
\[W_p(g)=G_p(g^{-1}\cdot S_g).\]
It is easy to check $\RR^s_g =W_{1}(g)$.
The $p$-Weitzenb\" ock curvature is as much interest as the scalar curvture,
at least
for \lcf manifold. For example, from the Weitzenb\" ock formula
one can easily show that  a  \lcf manifold with
positive $p$-Weitzenb\" ock curvature has vanishing
$q$-cohomology group for $n/2-p\le q\le n/2+p$ (for $p\le n/2$).

We may also ask whether we can find a conformal metric
with constant $W_p$. The corresponding equation is
\begin{equation}\label{eq3}
G_p\left(\uuu\right)=constant\cdot e^{-2u}.
\end{equation}
Equation (\ref{eq3}) is also uniformly elliptic. It is easy to see
that (\ref{eq3}) is concave for $p\le \frac n2$. From the
transformation formula (\ref{transf}), one can check that Equation
(\ref{eq3}) is also conformally invariant. Hence, in general there
is no compactness for equation (\ref{eq3}).

We also define for $p\ge n/2$ a function $G_p:\R^n\to \R$
by
\[G_p(A)=G_p(\L)=\max\{ (n-p)\sum _{k\le p}\l_{i_k}+p\sum _{k>p}\l_{i_k}\},\]
where $\max$ is over all permutations of $1,2,\cdots, n$. We still
define $W_p(g)=G_p(g^{-1}\cdot S_g)$. It is clear that ${\cal
R}^l_g=W_{n-1}(g)$. We can also consider equation (\ref{eq3}) for
$p\ge n/2$. For $p\ge n/2$ equation (\ref{eq3}) is still uniformly
elliptic, but it is convex. Since $G_p$ is homogeneous of order
$1$, we may rewrite equation (\ref{eq4.0}) as  the following
equivalent equation by setting $v=e^u$
\begin{equation}\label{neq4}
G_p(\n^2 v+S_{g_0}v)=p(n-p)\frac {|\n v|^2}v+\frac f v.
\end{equation}

The function $G_p$ is only Lipschitz, when $p\not =n/2$. (Remark
that when $p=n/2$, equation (\ref{eq3}) is equivalent to the
Yamabe equation.) One can find a sequence of smooth functions
$\{F_k\}$ such that $F_k$ uniformly converges to $G_p$
 in any compact domain of $\R^n$ and homogeneous $1$ outside the
 unit ball in $\R^n$, i.e., $F_k(x)=|x|F(\frac{x}{|x|})$ for
 $|x|\ge 1$. Furthermore, for $p<n/2$ ($p>n/2$) $F_k$ is concave
 (convex).

One may consider a more general class of conformal equations. Let
${\cal S}$ be the space of symmetric $2$-tensors on $M$. Let
$F:{\cal S} \to \R$ a real continuous function. We consider the
following general equation
\begin{equation}\label{geq}
F\left(e^{-2u}(\uuu)\right)=f(x),\end{equation} for some function
$f:M\to \R$. $F$ is  {\it uniformly elliptic} with constants
$\l_0$ and $\L_0$ if there exists two positive constants $\l_0$
and $\L_0$ such that for any $W\in {\cal S}$
\[\l_0  \|N\| \le F(W+N)-F(W)\le \L_0\|N\| \quad \forall N\ge 0,\]
here by $N\ge 0$ means that $N$ is nonnegative definite and
$\|N\|=\sup_{|v|=1}|Nv|$. If $F$ is uniformly elliptic, we call
equation (\ref{geq}) a uniformly elliptic fully nonlinear equation
with ellipticity constants $\l_0$ and $\L_0$. There are many
typical uniformly elliptic fully nonlinear equations. Our equation
(\ref{eq3}) is similar to the Pucci equation, see \cite{CC}. Let
${\cal M}^\pm$ be the Pucci's extremal operators, namely for two
given constant $0<\l_0<\L_0$ and $W\in {\cal S}$
\[\begin{array}{rcl}
 {\cal M}^-(W) &=&\ds\vs \l_0\sum_{e_i>0}e_i+\L_0\sum_{e_i<0}e_i,\\
 {\cal M}^+(W) &=&\ds \L_0\sum_{e_i>0}e_i+\l_0\sum_{e_i<0}e_i,\end{array}\]
where $e_i=e_i(W)$ are the eigenvalues of $W$. One can also consider
\begin{equation}\label{Pucci}
\begin{array}{rcl}
\ds {{\mathcal M}}^-(e^{-2u}(\uuu))&=&\ds\vs 1,\\
\ds {{\cal M}}^+(e^{-2u}(\uuu))&=&\ds 1.
\end{array}\end{equation}
${\cal M}^-$ is concave, while ${\cal M}^+$ is convex.

We now deduce $C^{2,\alpha}$ estimates for equation (\ref{eq4.0})
from the work of Caffarelli \cite{CC} assuming the gradient bound
and a lower bound of $u$. The crucial part in Caffarelli's
perturbation theory for uniformly elliptic fully nonlinear
equation $F(\nabla^2 u, x)=f(x)$ is the $C^{1,1}$ interior
estimates for $F(\nabla^2 u, x_0)=constant$. He obtained such
fundamental estimates for concave or convex operator $F$ (note
that here concavity and convexity of $F$ can always be switched by
$\tilde F(\L)=-F(-\L)$). Though Caffarelli proved these estimates
in \cite{CC} for equations with flat metric, his arguments work
under general Riemannian metrics. And the generalization of
Caffarelli's estimates by L. Wang \cite{LW} to uniformly elliptic
equations of form $F(\nabla^2u, \nabla u, u, x)$ gives the
following $C^{2,\alpha}$ estimates for equations of type
(\ref{eq3}).

\begin{thm}[$C^{2,\alpha}$-estimates] \label{thmC2} Suppose $F$ is a uniformly elliptic
concave operator with elliptic constants $\lambda_0, \Lambda_0$.
Let $B_1$ be a unit disk in a compact Riemannian manifold $M$ and
$f, h \in C^1(B_1)$. Suppose $g=e^{-2u}g_0$ with $|\nabla
u|_{B_1}\le A$ is a solution of equation
\begin{equation}\label{ng1} F\left(\uuu
\right)=e^{-2u}f(x), \quad x\in B_1,\end{equation} then there
exist $\alpha>0$ and $C>0$ depending only on $\lambda_0,
\Lambda_0, A, \|u\|_{C^0(B_1)}$ and $g_0$ such that
\begin{equation}\label{eqCC5}
\| u\|_{C^{2,\alpha}(B_{\frac12})} \le C.
\end{equation}
\end{thm}

In fact, we may directly apply Caffarelli's estimates \cite{CC} to
obtain $C^{2,\alpha}$ estimates for equations of type
(\ref{neq4}). Let $B_1$ be a unit disk in a compact Riemannian
manifold $M$ and $f, h \in C^1(B_1)$. Suppose
$g=e^{-2u}g_0=v^{-2}g_0$ with $|\frac {\nabla v}{v}|_{B_1}\le A$
is a solution of equation
\begin{equation}\label{nng1} F(\n^2
v(x)+v(x)S_{g_0}(x))=h(x)\frac{|\n v(x)|^2}{v(x)}+\frac
{f(x)}{v(x)}, \quad x\in B_1.\end{equation} Since $F$ is concave,
by Theorem 6.6 in \cite{CCB}, the equation
\[F(\nabla^2 v+S_{g_0}(x_0)v)=constant,\]
has $C^{1,1}$ interior estimates for any $x_0\in B_1$. It follows from
Theorem 7.1 in \cite{CCB} that equation (\ref{nng1}) has interior
$W^{2,p}$ estimate for any $n<p<\infty$ since $\frac{|\nabla
v|}{v} \le A$. This in turn gives $C^{1,\beta}$ a priori bound for
the solution $v$ of equation (\ref{nng1}) for all $0<\beta<1$. Finally
estimate (\ref{eqCC5}) for $v=e^u$ follows from Theorem 8.1 in
\cite{CCB} since the right hand side of equation (\ref{nng1}) is a
$C^{\beta}$ function now.

\section{Local gradient estimates}

In this section, we establish a local gradient estimates for
solutions $g=e^{-2u}g_0$ of equation (\ref{eq4.0}) in the case
$p\le n/2$, depending only on lower bound of $u$.

\begin{thm}[Local gradient estimates] \label{thm1}Suppose $F$ is concave and uniformly
elliptic with ellipticity constants $\lambda_0, \Lambda_0$. Let
$B_1$ be a unit disk in a compact Riemannian manifold $M$ and $u$
a $C^{2}$ solution of the following equation
\begin{equation}\label{eq4}F\left(\uuu\right)=e^{-2u}f(x), \quad x \in B_1 \end{equation}
for a $C^1$ function $f:B_1\to \R$. Then there is a constant $C>0$
depending only on $\lambda_0, \Lambda_0, g_0$ such that
\begin{equation}\label{eq5}
|\n u|^2(x) \le C (1+ \|f\|_{C^1(B_1)}e^{-2\inf_{B_1} u}), \quad
\hbox{ for any } x\in B_{1/2}.
\end{equation}
\end{thm}

Combining Theorem \ref{thm1} and Theorem \ref{thmC2}, we deduce
the following.
\begin{cor}[Local $C^{2,\alpha}$-estimates] \label{nthmC2}
Let $B_1$ be a unit disk in a compact Riemannian manifold $M$ and
$f\in C^1(B_1)$. Suppose $u$ is a solution of equation
(\ref{eq4}), then there is a constant $C>0$ depending only on
$\lambda_0, \Lambda_0, g_0, \|f\|_{C^1{B_1}}, \inf_{B_1}u,$ such
that
\begin{equation}\label{eqC5}
\| u\|_{C^{2,\alpha}(B_{\frac12})} \le C.
\end{equation}
In particular, (\ref{eqC5}) is true for any solution of
(\ref{eq4.0}) when $p\le \frac n2$.
\end{cor}

We note that Corollary \ref{nthmC2} guarantees that we can use the
so-called blow-up analysis to study the fully nonlinear equation
(\ref{eq4.0}). We remark that the local gradient estimates
(\ref{eq5}) does not true for $p=n-1$. The local gradient
estimates were first obtained for $\s_k$ conformal invariant
equations in \cite{GW1} and for $\frac{\s_k}{\s_l}$ conformal
invariant equations in \cite{GW3}, though in general fully
nonlinear equations have no local estimates. The operator $G_p$ we
are considering here is only Lipschitz, which we will deal with by
a smoothing argument. Actually we can prove the local gradient
estimates for a more general class of uniformly elliptic fully
nonlinear conformal equations.

We first prove the local gradient estimates for $C^2$ uniformly
elliptic operator $F$. Let $F: \R^n\to \R$ is a $C^2$ symmetric
function and consider the following equation
\begin{equation}\label{eq3.1}
F\left(\uuu\right)=\tilde f,\end{equation} for some $C^1$ function
$\tilde f$. We denote the left hand side  of (\ref{eq3.1}) by
$F(W)$ and set
\[F^{ij}=\frac{\partial F}{\partial w_{ij}},\]
where $w_{ij}$ is the entry of the matrix $W$. As mentioned above,
$F(W) =F(\L)$, where $\L$ is the set of eigenvalues of $W$.

\begin{pro}\label{lem1}Let $B_1$ be a unit disk in a compact Riemannian
manifold $M$ and $u$ a $C^{3}$ solution of (\ref{eq3.1}) for a
$C^1$ function $\tilde f:B_1\to \R$. Let $F:{\cal S}\to \R$ be a
$C^2$ function satisfying
\begin{itemize}
\item[(1)]  $F$ is an uniformly elliptic with ellipticity
constants $\l_0$ and  $\L_0$ \item[(2)] $F$ is concave.
\end{itemize}
Then for any $\rho \in C^2_0(B_1)$ with $0\le \rho(x)\le 1$, there
is a constant $C>0$ depending only on $\l_0$, $\L_0$,
$\|\rho\|_C^2(B_1)\|$ and $g_0$ such that
\begin{equation}\label{eq3.3}
\max_{B_1}\{\rho(x)|\n u|^2(x)\} \le C (1+\max_{B_1}\{\rho(x)|\n
\tilde f(x)|\}).\end{equation}
\end{pro}

\

\noindent{\it Proof of Proposition \ref{lem1}.} The Proof follows
closely the argument given in \cite{GW1} and \cite{GW2}. As in
\cite{GW1}, we first reduce the proof of the Lemma to the
following  claim.
 \medskip

 \noindent{\bf Claim.} {\it There is a constant $A_0$ depending,
such that
\begin{equation}\label{B1}
\sum_{i,j}F^{ii}\tilde u^2_{ij} \ge A_0^{-\frac 58}\sum_i
F^{ii}|\n u|^4,
 \end{equation}}
\noindent where $\tilde u_{ij}=u_{ij}+S_{ij}$.

\

For convenience of the reader, we sketch the reduction. Let $\rho$
be a test function $\rho\in C^2_0(B_1)$. We may assume
\begin{equation}\label{8}
|\n \rho (x)| \le 2 b_0 \rho^{1/2}(x) \quad \hbox{ and }\quad
|\n^2 \rho|  \le   b_0,    \hbox{ in } B_1, \end{equation} for
$b_0>1$. Set $H=\rho |\n u|^2$. Our aim is to bound $\max_{B_1}H$.
Let $x_0\in B_1$ be a maximum point of $H$ and assume that $W$ is
a diagonal matrix at the point $x_0$ by choosing a suitable normal
coordinates around $x_0$. Set $\l_i=w_{ii}$ and
$\L=(\l_1,\l_2,\cdots,\l_n)$. Since $W$ is diagonal at $x_0$, we
have at $x_0$
\begin{equation}\label{add0}
    w_{ii}  =  u_{ii}+u^2_i-\frac12 |\n u|^2+S_{ii}, \quad
    u_{ij}  = -u_iu_j-S_{ij}, \quad \forall i\neq j,
\end{equation}
where $S_{ij}$ are entries of $S_{g_0}$. We may assume that
\[H(x_0)\ge b_0^2A_0^2,\]
for some large, but fixed constant $A_0>0$
which will be fixed later. We may also assume
that
\begin{equation}\label{add1.1} |S_{g_0}| (x_0)\le A_0^{-1}|\n
u|^2(x_0).\end{equation} Otherwise, we are done. The fact
that the derivatives of $H$ at
$x_0$ vanish  imply
\begin{equation}\label{add2.2}
\left|\sum_{l=1}^nu_{il}u_l\right|(x_0)\le \frac{|\n
u|^3}{A_0}(x_0) \quad \hbox{ for any }i.\end{equation}
Applying the maximum principle to $H$, we have
\begin{equation}\label{add100}
0 \ge \ds\vs
 F^{ij}H_{ij}=F^{ij}
\left\{\left(-2\frac{\rho_i\rho_j}\rho+\rho_{ij}\right)|\n u|^2
+2\rho u_{lij}u_l+2\rho u_{il}u_{jl}\right\}.\end{equation} The
first term in the left hand side of (\ref{add100}) is bounded from
below by $-10nb_0\Lambda_0|\n u|^2$. By using equation
(\ref{eq3.1}) and inequality (\ref{add2.2}),
 the second term can be bounded by
\begin{equation}\label{14}\ba{rcl}
\ds\vs \sum _{i,j,l}F^{ij}u_{ijl}u_l & \ge & \ds
 \sum_lF_lu_l
-2\sum_{i,l}F^{ii}u_{il}u_lu_i+\sum_{i,l}
 F^{ii}u_{jl}u_ju_l
 -C|\n u |^2 \sum_{i}F^{ii}\\
&\ge& \ds -|\tilde f|^2 -2n  \sum_{i} F^{ii} \frac {\ds |\n
u|^4}{\ds A_0} -C|\n u |^2 , \ea\end{equation} where $C>0$ depends
only on $g_0$ and $\Lambda_0$. See also (2.20) in \cite{GW1}. It
is easy to see that the third term is bounded by the Claim. Hence
if the Claim is true, from (\ref{add100}) we have
\begin{equation}\label{add101}
0\ge -C|\n u|^2-\rho |\n \tilde f|^2-\rho \sum_i F^{ii} \frac{|\n
u|^4}{A_0}+\rho A_0^{-\frac 5 8} \sum_iF^{ii}|\n
u|^4.\end{equation} Multiplying (\ref{add101}) by $\rho$, we have
\[0\le \sum_iF^{ii}(A_0^{-\frac 58}-A_0^{-1})
H^2 -CH-\rho^2 |\n \tilde f|^2,\] from which we have
(\ref{eq3.3}).

Now we prove the Claim. By (\ref{add0}), we have
\begin{equation}\label{19}\ba{rcl}
\vs \ds \sum_{i,l}F^{ii}\tilde u_{il}^2
& =
& \ds \vs \sum_{i}F^{ii}\tilde u_{ii}^2+ \sum_{i\not = l}F^{ii}u_i^2 u^2_l
\\
& =& \ds \vs \sum_{i}F^{ii}\left\{\tilde u_{ii}^2 +u_i^2(|\n u|^2-u_i^2)
\right\}\\
&= &  \ds  \sum_{i}F^{ii} (w^2_{ii}-2u^2_iw_{ii}+w_{ii}|\n u|^2+
\frac {|\n u|^4}4).\\
\ea\end{equation}
Set $\d_0= A_0^{-1/4}<0.1$. We divide the set $I=\{1,2,\cdots,
n\}$ as in \cite{GW1} into two parts:
\[I_1=\{i\in I\,|\,u_i^2 \ge \d_0 |\n u|^2\} \quad \hbox{ and } \quad
I_2=\{i\in I\,|\,u_i^2 < \d_0 |\n u|^2\}.\] It is clear that $I_1$ is
non-empty.

\

{\noindent \bf Case 1.}  There is $j_0$ satisfying
\begin{equation}\label{B4}
\tilde u_{jj}^2 \le \d_0^2 |\n u|^4 \quad \hbox{ and } \quad
u^2_{j} <\d_0|\n u|^2.\end{equation}

\

We may assume that $j_0=n$.
We have $|w_{nn}+\frac {|\n u|^2} 2|= |\tilde u_{nn}+ u_n^2|<2\d_0
|\n u|^2$ by (\ref{B4}). From (\ref{add2.2}) and (\ref{add0}), we
have
\[\begin{array}{rcl} \left|w_{ii}-\frac
{|\n u|^2} 2\right|&=& \left|u_{ii}+u_i^2-|\n u|^2+S_{ii} \right|
\le 3\d_0^2 |\n u|^2,\end{array}\] for any $i\in I_1$.

Using these estimates, we  repeat the derivation of equation
(2.38) in \cite{GW1} to obtain
\begin{equation}\label{B3}
\ds\vs \sum_{i,l} F^{ii} \tilde u^2_{il}\ge
\tilde F^1\frac {|\n u|^4} 4-\tilde F^1 |\n u|^4+
F^{nn}\frac {|\n u|^4} 4+(1-32 \d^2_0) \frac{|\n u|^4} 4\sum_{i} F^{ii},
\end{equation}
where
$\tilde F^1:=\max_{i\in I_1} F^{ii}$.
Recall that $I_1$ is necessarily non-empty. We  assume $1\in I_1$ with
 $F^{11}=\tilde F^1$. The concavity of $F$ implies that
 \begin{equation}\label{cvx} F^{nn} \ge F^{11},\end{equation}
 for $w_{11} >w_{nn}$. Hence, from (\ref{B3}) we have
 \[\begin{array}{rcl}
\ds\vs \sum_{i,l} F^{ii} \tilde u^2_{il}&\ge&\ds
-16\d_0^2|\n u|^4+(1-32 \d^2_0) \frac{|\n u|^4} 4\sum_{i=2}^{n-1}
 F^{ii}\\
 &\ge &\ds \vs \left(1-32 \d^2_0-64 \d_0^2\frac {\L_0}{(n-2)\l_0}\right) \frac{|\n u|^4} 4\sum_{i=2}^{n-1}
 F^{ii}\\
&\ge & \ds \left(1-32 \d^2_0-64 \d_0^2\frac {\L_0}{(n-2)\l_0}\right) \frac {(n-2)\l_0}{n\L_0}
 \frac{|\n u|^4} 4\sum_{i=2}^{n-1}
 F^{ii},
\end{array}\]

\

{\noindent \bf Case 2.} There is no $j\in I$ satisfying (\ref{B4}).

\medskip

For this case, the proof is the same as in \cite{GW2}.
We repeat it here for completeness.

We may assume that there is $i_0$ such that $\tilde
u_{i_0i_0}^2\le \d_0^2|\n u|^4$, otherwise the claim is
automatically true. Assume $i_0=1$. As in Case 4 in \cite{GW1}, we
have $u_1^2\ge (1-2\d_0)|\n u|^2$ and $\tilde u^2_{jj} + u_j^2(|\n
u|^2-u_j^2)\ge \d^2_0 |\n u|^4$ for $j>1$. From equation
(\ref{19}), we have
\begin{eqnarray*}
\sum_{i,l}F^{ii}\tilde u_{il}^2 & =& \sum_{i}F^{ii}\left\{\tilde u_{ii}^2 +u_i^2(|\n u|^2-u_i^2)
\right\}\\
& \ge &  \sum_{i\ge 2}F^{ii}(\tilde u^2_{ii}+u_i^2(|\n u|^2-u_i^2))\\
&\ge& \d^2_{0} |\n u|^4\sum_{i\ge 2} F^{ii} \ge C
\d^2_{0} |\n u|^4\sum_{i\ge 1} F^{ii}.
 \end{eqnarray*}
 The latter inequality follows from the uniformly ellipticity of
$F$. This finishes the proof the Claim and hence the Proposition.
\qed

\

We have a direct corollary.
\begin{cor}\label{thmn100}Suppose $F$ is a $C^2$ concave and uniformly
elliptic operator with ellipticity constants $\lambda_0,
\Lambda_0$. Let $B_1$ be a unit disk in a compact Riemannian
manifold $M$ and $u$ a $C^{2}$ solution of equation
\begin{equation}\label{eqn104}F\left(\uuu\right)=e^{-2u}f(x), \quad x \in B_1 \end{equation}
for a $C^1$ function $f:B_1\to \R$. Then there is a constant $C>0$
depending only on $\lambda_0, \Lambda_0, g_0$ such that
\begin{equation}\label{eqn105}
|\n u|^2(x) \le C (1+ \|f\|_{C^1(B_1)}e^{-2\inf_{B_1} u}), \quad
\hbox{ for any } x\in B_{1/2}.
\end{equation}
\end{cor}

\pr We pick $\rho \in C^2_0(B_1)$ such that $\rho(x)=1, \forall
x\in B_{\frac 12}$ and $0\le \rho(x)\le 1, \forall x\in B_1$.
(\ref{eqn105}) follows directly from (\ref{eq3.3}) with $\tilde
f=e^{-2u}f$. \qed

In what follows in the next sections, we will only need Corollary
\ref{thmn100} as we will work on smooth operator $F$. We note
that estimates (\ref{eqn105}) and (\ref{eqC5}) are independent of
the smoothness of $F$. Theorem \ref{thm1} can also be proved by
certain appropriate approximations.

\medskip

\noindent {\it A sketch proof of Theorem \ref{thm1}.} Since $u\in
C^2$, $u$ is in fact $C^{2,\alpha}$ by the Evans-Krylov theorem.
We may find two sequences of smooth functions $\{u_k\}$ and
$\{f_k\}$, such that $u_k \to u$ in $C^{2,\alpha}(\bar B_1)$,
$f_k\to f$ in $C^{0,1}(\bar B_1)$, and \[F\left(\n
^2u_k+du_k\otimes du_k-\frac {|\n u_k|^2}2g_0+S_{g_0}\right)\ge
e^{-2u_k}f_k.\] We now construct a sequence of smooth concave
$F_k:{\cal S}\to \R$ such that $F_k$ converges to $F$ uniformly in
compacts of $\cal S$ and $F_k$ is uniformly elliptic with
ellipticity constants $\frac{\l_0}2$ and $2\L_0$. We may assume
\begin{eqnarray*}
 F_k(\l)\ge F(\l), \quad
 \forall |\l|\le \sup_{B_1}
\left|\uuu\right|+1.\end{eqnarray*} By the symmetry of $F_k$, $\frac{\partial
F_k(1,\cdots,1)}{\partial \l_i}=\frac{\partial
F_k(1,\cdots,1)}{\partial \l_j}$ for all $i,j$. Set
$\frac{\partial F_k(1,\cdots,1)}{\partial \l_j}=A$. Let $R_0(x)$
be the scalar curvature of $g_0$, and let $\tilde u$ be the
solution of $\Delta \tilde
u_k=\frac{f_ke^{-2u_k}}A-R_0+n-\frac{F_k(1,\cdots,1)}{A}$ in $B_1$
with $\tilde u_k=u_k$ on $\partial B_1$. We consider the following
Dirichlet problem
\begin{equation}\label{Diri}\begin{array}{rcll}
\ds\vs F_k\left(\n ^2v_k+dv_k\otimes dv_k-\frac {|\n
v_k|^2}2g_0+S_{g_0}\right) &=&\ds f_ke^{-2u_k} , &\quad
\hbox{ in } B_1,\\
v_k&=&u_k, &\quad \hbox{ on }\partial B_1.
\end{array}\end{equation}
By the concavity of $F_k$,
\[F_k(\l)\le A\sigma_1(\l) -nA+F_k(1,\cdots,1).\]
We have
\begin{eqnarray*}  A \left(\Delta v_k+R_0-n+\frac{F_k(1,\cdots,1)}{A}\right)
&\ge & A\left(\Delta v_k -\frac{n-2}{2}|\nabla
v_k|^2+R_0-n+\frac{F_k(1,\cdots,1)}{A}\right) \\
&\ge & f_ke^{-2u_k} \\
&\ge & F_k\left(\n ^2v_k+dv_k\otimes dv_k-\frac {|\n
v_k|^2}2g_0+S_{g_0}\right).\end{eqnarray*} In turn, we have $v_k\le \tilde
u_k $ in $B_1$. On the other hand, we have
\begin{eqnarray*}  F_k\left(\n
^2u_k+du_k\otimes du_k-\frac {|\n u_k|^2}2g_0+S_{g_0}\right) &\ge
& F\left(\n
^2u_k+du_k\otimes du_k-\frac {|\n u_k|^2}2g_0+S_{g_0}\right)\\
&=&F_k\left(\n ^2v_k+dv_k\otimes dv_k-\frac {|\n
v_k|^2}2g_0+S_{g_0}\right).\end{eqnarray*} This gives $v_k\ge u_k$
in $B_1$. From this, we obtain a $C^0$ bound of $v_k$ and a bound of
$|\nabla v_k|$ at the boundary $\partial B_1$. Using the same
proof of Proposition \ref{lem1}, we can obtain a bound of $|\nabla
v_k|$ on $\bar B_1$ (simply let $H(x)=|\nabla v_k|^2$ and estimate at
the maximum point if it is not on the boundary). At this end, we
have a uniform $C^1$ bound of $v_k$. The standard barrier
construction $\omega_{\pm}$ similar to the one in Step 3 in
Chapter 9 of \cite{CCB} (page 91), with the modified operator
$\tilde F(\omega_{\pm})=F_k(\n ^2 \omega_{\pm}+dv_k\otimes
dv_k-\frac {|\n v_k|^2}2g_0+S_{g_0})$, will give a $C^2$ bound near
boundary. The global $C^2$ estimate follows easily along the lines
of proof in Proposition 3.1 in \cite{GW1} (see also proof of Lemma
3 in \cite{GLW}). Higher regularity estimates follow from the
Krylov Theorem \cite{Krylov}. We note that $C^2$ and higher
regularity bounds of $v_k$ may depend on higher smoothness
assumptions on $u_k$ and $f_k$, but the interior $C^{2,\alpha}$
estimates of $v_k$ depend only on $\|u_k\|_{C^{2,\alpha}(\bar
B_1)}$ and $\|f_k\|_{C^1(\bar B_1)}$ by Theorem \ref{thmC2}. In
any case, we can establish the existence of the Dirichlet problem
(\ref{Diri}) by using the method of continuity for equation
\[F_t\left(\n ^2v_k+dv_k\otimes dv_k-\frac {|\n
v_k|^2}2g_0+S_{g_0}\right)= tf_ke^{-2u_k}+(1-t)f_k^*,\] where
$F_t(\l)=tF_k(\l)+(1-t)\sigma_1 (\l)$ and $f_k^*=\sigma_1\left(\n
^2u_k+du_k\otimes du_k-\frac {|\n u_k|^2}2g_0+S_{g_0}\right)$. Now
by Proposition \ref{lem1}, we have for any $\rho \in C^2_0(B_1)$
with $0\le \rho \le 1$, there exists $C$ independent of $k$ and
$u_k$ such that
\begin{equation}\label{a1}
\max_{B_1}\{|\rho(x)\n v_k|^2(x)\}  \le
C(1+(\max_{B_1}\{|\rho(x)f_k(x)||\n
u_k(x)|\}+\max_{B_1}\{\rho(x)|\n f_k(x)|\})e^{-2u_k(x)}).
\end{equation}
 From our estimates, $\|v_k\|_{C^1(\bar B_1)}\le C$, where $C$ is  a constant
independent of $k$, since $\|u_k\|_{C^{2,\alpha}(\bar B_1)}$
and $\|f_k\|_{C^1(\bar B_1)}$ are uniformly bounded. By Theorem
\ref{thmC2}, $v_k \to v_0$ (after passing a subsequence) in
$C^{2,\alpha}(B_1)$, $v_k$ converges to $v_0=u$ by the uniqueness.
Therefore, there is a constant $C>0$ depending only on $\l_0$,
$\L_0$, and the geometry of $B_1$ such that
\[\max_{B_{\frac 12}}|\n u|^2 (x)\le C(1+\|f\|_{C^1(B_1)} e^{-2\inf_{B_1}u(x)}).\]
This finishes the proof of Theorem \ref{thm1}. \qed

\

A direct consequence of Theorem \ref{thm1} is the following
\begin{cor}\label{coro1} Let $B_1$ be a unit disk in a Riemannian manifold
$(M, g_0)$ and $p<n/2$. There exists a small constant $\e_0>0$ depending only
 on $(B_1,g_0)$
such that for any sequence of solutions $u_i$ of (\ref{eq4.0}) in
$B_1$ with
\[\int_{B_1}e^{-nu_i} dvol(g_0)\le \e_0,\]
either
\begin{itemize}
\item[(1)] there is a subsequence $u_{i_l}$ uniformly converging
to $+\infty$ in any compact subset of $B_1$, or
\item[(2)] there is a subsequence $u_{i_l}$ converging strongly in $C^{2,\a}_{loc}
(B_1).$
\end{itemize}\end{cor}
\pr It follows the same lines of proof in \cite{Sch} or
\cite{GW1}.\qed

\section{negative curvature case}
In this section, we discuss the negative curvature case, where the
geometry is rich. By \cite{GY} ($n=3$) and \cite{Lo} (general
dimension $n\ge 3$), every higher dimensional manifold has a
metric with negative Ricci tensor. It is clear that such a metric
also has negative $W_p$ for any $1\le p\le n-1$. Hence every
higher dimensional manifold has a metric with negative $W_p$. The
conformal deformation will yield interesting geometric information
about the extremal metrics in a given conformal class.

\medskip

\noindent{\it  Proof of Theorem \ref{thm2}.} First we take the
sequence of smooth $F_k$ considered  in the previous section. For
each $F_k$ we consider the following equation
\begin{equation}\label{eq4.1}
F_k \left((e^{2u}(\uuu)\right)=-1.\end{equation} For large $k$,
from the condition of the Theorem, we have
\[F_k(S_{g_0})(x)<0, \quad \forall x\in M.\]
We first prove the existence of solutions to equation (\ref{eq4.1}).
Here we use the method of continuity. Let us consider
the following equation
\begin{equation}\label{eq4.2}
F_t(u):=F_k \left(
e^{2u}(\uuu)\right)+t-(1-t)F_k(S_{g_0})=0\end{equation} and define
$J=\{t\to [0,1]\,|\, \hbox{ (\ref{eq4.2}) has  a solution for }
t\}.$ It is clear that $0\in J$. First, we prove the openness of
$J$. Let $t_0\in J$. By the maximum principle, we know that there
is only one solution $u$ of (\ref{eq4.2}) for $t=t_0$. Let $L$ be
its linearization. We want to show that $L$ is invertible. By the
maximum principle again, we know that the kernel of $L$ is
trivial. Note that $L$ might be not self-adjoint. To show the
invertibility of $L$, we  need to show that the cokernel of $L$ is
also trivial. However, one can readily check that the Fredholm
index of $L$ is zero, and hence the cokernel of $L$ is trivial.
Now the openness follows from the implicit function theorem.

Then we show the closeness. Let $x_0$ and $x_1$ be  the minimum
and maximum of $u$ respectively. By the maximum principle, we have
\begin{equation}\label{eqC0}
e^{2u(x_0)}\ge \frac{t-(1-t)F_k(S_{g_0})(x_0)}{-F_k(S_{g_0})(x_0)}\end{equation}
and
\begin{equation}\label{eqC00}
e^{2u(x_1)}\le \frac{t-(1-t)F_k(S_{g_0})(x_1)}{-F_k(S_{g_0})(x_1)}.
\end{equation}
Hence, we have $C^0$ bound of $u$ independent of $t$.

By a  global estimates proven in Proposition \ref{pro4.1} below
and Theorem \ref{thmC2}, we have the closeness. Hence we have a
solution $u_k$ of (\ref{eq4.1}) with the bound
\begin{equation}\label{eqC1}
\frac{1}{-\min F_k(S_{g_0})(x_0)} \le e^{2u_k} \le \frac{1}{-\max
F_k(S_{g_0})(x_0)}.\end{equation} In viewing of (\ref{eqC1}), we
use again the global estimates and Theorem \ref{thmC2} to obtain a
$C^{2,\a}$ uniform bound of $u_k$ for some $\a>0$. Hence $u_k$
converges (by taking a subsequence) to $u$. It is clear that
$e^{-2u}g_0$ satisfies (\ref{eq4.0}).\qed

\begin{pro}\label{pro4.1} Let $u$ be a solution of
\begin{equation}\label{eqG1} F( \uuu)=fe^{-2u}\end{equation}
with $C^1$ function $f$. Suppose that  $F$ is uniformly elliptic
and is homogeneous of degree $1$. Assume that $ u$ has $C^0$
bound. Then $u$ has a $C^1$ bound and a $C^{2,\a}$ bound.
\end{pro}

\pr Since we already have full $C^0$ bound, the proposition can be
proved using standard Pogorelov type of trick, for example, as in
\cite{GV2}. Let $v=e^u$ and consider the following equivalent form
of (\ref{eqG1})
\begin{equation}\label{eqG2}
F\left(\frac{\n^2 v}v -\frac{|\n
v|^2}{v^2}g_0+S_{g_0}\right)=fv^{-2}.\end{equation} Without loss
of generality, we may assume that $v\le 1$. Set $H=e^{2\phi(v)}|\n
v|^2$. Here $\phi$ will be fixed later. Let $x_0$ be a maximum
point of $H$. At $x_0$, we have
\begin{equation}\label{eqpro4.1.1}
\sum_k (2v_kv_{ki}+2\phi'(v)  v_i |\n v|^2)=0, \quad \hbox{ for
any } i.\end{equation} Without loss of generality, we may assume
that $v_1=|\n v|$, $v_i=0$ for any other $i$ and that $v_{ij}$ is
diagonal at $x_0$. Hence, (\ref{eqpro4.1.1}) is equivalent to
$v_{11}=-\phi'(v)|\n v|^2$. Set $w_{ij}=v_{ij}-\frac{|\n
v|^2}{2v}\d_{ij}+vS_{ij}$. The maximum principle, together with
(\ref{eqpro4.1.1}), implies
\begin{equation}\label{eqpro4.1.2}
0\ge \ds \sum F^{ij} (v_{kj}v_{ki}+ v_kv_{kij}+\phi'(v) v_{ij}|\n
v|^2+v_iv_j\phi''(v)|\n v|^2
+2\phi'(v)v_{i}v_{ki}v_k)\end{equation} By (\ref{eqpro4.1.1}), it
is easy to check that
\[ \sum F^{ij} v_{ki}v_{kj} \ge F^{11}
v_{11}v_{11}=F^{11} (\phi'(v))^2\frac{|\n v|^4}{v},\]
\[ \begin{array}
{rcl} \ds  \sum F^{ij} v_k v_{kij} & =&\ds \vs \sum F^{ij}  v_k
v_{ijk} +
\sum F^{ij}v_kv_mR^m_{ijk}\\
&\ge&\ds\vs \sum F^{ij}v_k (w_{ij}+\frac {|\n v|^2} {2v}\d_{ij}+vS_{ij})_k
-C \max v|\n v|^2\\
&\ge&\ds\vs 2v_k(f v^{-1})_k+\sum F^{ii}v_1(-\frac{v_1^3}{2v^2}-\phi'(v)
\frac{|\n v|^3}{v})-C \max v|\n v|^2\\
&\ge&\ds\vs 2 f_k v_k-2v^{-1}|\n v|^2-(\frac 1 {2v}+\phi'(v)) \sum
F^{ii}\frac{|\n v|^4} v-C \max v|\n v|^2,\end{array}\]

\[ \begin{array}
{rcl} \ds \phi'(v)\sum F^{ij} v_{ij}|\n v|^2 &=&\ds \vs \phi'(v)
\sum F^{ij}(w_{ij}+\frac{|\n v|^2}{2v}\d_{ij}-vS_{ij})
|\n v|^2\\
&\ge&\ds\vs \phi'(v)fv^{-1}|\n v|^2+\phi'(v)\sum F^{ii}\frac{|\n
v|^4} {2v} -C \max v|\n v|^2, \end{array}\]

\[ \phi''(v)\sum F^{ij} v_iv_j|\n v|^2 =\phi''(v)
F^{11} |\n v|^4,\]

\[ 2\phi'(v)\sum F^{ij} v_jv_kv_{ki} = 2\phi'(v) F^{11}
v_1^2v_{11}=-2(\phi'(v))^2 F^{11}|\n v|^4.
\]
Here $C$ is a positive constant depending only on the Riemannian
curvature
of the background metric and it varies from line to line.
Therefore, we have
\begin{equation}\label{eqc1}
\begin{array}{rcl}
\ds\vs 0 &\ge& \ds F^{11}
(\phi''(v)-(\phi'(v))^2) |\n v|^4+\frac{1}{2v}\sum F^{ii}(-\frac 1  {v}+\phi(v))
|\n v|^4 -C\max v|\n v|^2\\
&&+\ds 2f_kv_k-2v^{-1}|\n v|^2+\phi'(v) f v^{-1}|\n v|^2.\end{array}\end{equation}
Choose
\[\phi=-\frac 12 \log t(ct-2),\]
for a large constant $c>0$ so that $c\min v>3$.
One can easily check that for any $t \in [\min v, \infty)$
\[\phi''(v)-(\phi'(v))^2=\frac 1{t^2(ct-2)^2}>0\]
and
\[
-(\frac 1 {t}+{\phi(t)}) =\frac 1{t(ct-2)}.\]
In view of (\ref{eqc1}), we have
\[0\ge \frac1{2(\max v)^2(c\max v-2)}|\n v|^4F^{ii}-c(f,C)|\n v|^2.\]
Now we have a global bound of $|\n v|$, which depend only on
$g_0$, $f$, $\min v$ and $\max v$. The $C^{2,\a}$ bound follows
from Theorem \ref{thmC2}.\qed


\

\begin{rem}\label{rm1} The condition $W_p(g_0)(x)<0, \forall x\in M$ in
Theorem \ref{thm2} can be weaken to $W_p(g_0)(x)\le 0, \forall
x\in M$ and $W_p(g_0)(x_0)< 0$ for some $x_0\in M$. In fact, under
the weaker condition, one may produce a metric $g\in [g_0]$ with
the stronger condition holds. This can be done using the short
time existence of the fully nonlinear flow \[u_t=G_p(\uuu)e^{2u},
\quad u|_{t=0}=0.\] The short time existence follows from standard
nonlinear parabolic theory, and the strict negativity of $W_p(g)$
(which is equal to $u_t$) follows from the strong maximum
principle.
\end{rem}
\begin{rem}It is of interest to characterize when the condition in
Theorem \ref{thm2} is true by some conformal geometric quantities.
The difficulty here is the lack of variational structure for this
type of equations. We note that when $p> \frac n2$, if the Yamabe
constant $Y([g_0])$ of $(M,g_0)$ is non-positive, then the
condition in Theorem \ref{thm2} is satisfied unless $(M,g_0)$ is
conformally equivalent to a Ricci flat manifold. This simple
observation follows from the solution of the Yamabe problem and
the fact that if the scalar curvature vanishing identically, for
$p> \frac n2$, $G_p\le 0$ and $G_p$ vanishes identically if and
only if the metric is Ricci flat.
\end{rem}

\begin{rem} It is also an interesting problem to consider the equation
$W_p(g)=-1$ on a complete, non-compact manifold. The arguments in
the proof of the existence of Dirichlet problem (\ref{Diri}) can
be extended to deal with a given boundary condition at the
infinity for equation $W_p(g)=-1$ on a complete non-compact
negatively curved manifold.
\end{rem}

\

As a direct consequence of Theorem \ref{thm2}, we have

\begin{cor}\label{corng} If there is $\tilde g\in [g_0]$ with
${\mathcal R}_{\tilde g}^s(x)< 0$ for all $x\in M$ then there is a unique $C^{2,\a}$
metric
$g^*\in [g_0]$   for some $\alpha>0$
such that ${\mathcal R}_{g^*}^s(x)=-1, \forall x\in M$. Similarly,
If there is $\tilde g\in [g_0]$ with ${\mathcal R}_{\tilde
g}^l(x)< 0$ for all $x\in M$, then there is a unique $g_*\in [g_0]$ and
$g_*\in C^{2,\alpha}(M)$ for some $\alpha>0$ such that ${\mathcal
R}_{g_*}^l(x)=-1, \forall x\in M$.\end{cor}

\

Corollary \ref{corng}  can be applied to consider minimal volumes
in conformal classes.
Set
\[{\mathcal C}^{-}=\{g\in [g_0] | {\mathcal R}_g^s(x)\ge -1, \forall
x\in M\},\]
\[{\mathcal C}_{-}=\{g\in [g_0] | {\mathcal R}_g^l(x)\ge -1, \forall
x\in M\},\]
\[\tilde{\mathcal C}^{-}=\{g\in [g_0] | {\mathcal R}_g^s(x)\le -1, \forall
x\in M\},\]
\[\tilde {\mathcal C}_{-}=\{g\in [g_0] | {\mathcal R}_g^l(x)\le -1, \forall
x\in M\}.\] Define
\[{\mathcal V}^s([g_0])=\inf_{g\in {\mathcal C}^{-}} vol(g),
\qquad {\mathcal V}^l([g_0])=\inf_{g\in {\mathcal C}_{-}}
vol(g),\]
\[{\mathbb V}^s([g_0])=\sup_{g\in \tilde {\mathcal C}^{-}} vol(g),
\qquad {\mathbb V}^l([g_0])=\sup_{g\in \tilde{\mathcal C}_{-}}
vol(g),\]
where $vol(g)$ is the volume of $g$.
It is clear that the definitions  given here

\begin{lem}\label{ncomp} Let $g_1,g,g_2\in [g_0]$. If
${\mathcal R}_g^s(x)<0,$ and  ${\mathcal R}_{g_1}^s(x)\le
{\mathcal R}_g^s(x)\le {\mathcal R}_g^s(x), \forall x\in M$, then
$vol(g_1)\le vol(g)\le vol(g_2)$, any one of the equalities
holds if and only if the metric is equal to $g$. Similarly, if
${\mathcal R}_g^l(x)<0,$ and  ${\mathcal R}_{g_1}^l(x)\le
{\mathcal R}_g^l(x)\le {\mathcal R}_g^l(x), \forall x\in M$, then
$vol(g_1)\le vol(g)\le vol({g_2})$, any one of the equalities
holds if and only if the metric is equal to $g$.\end{lem}

The Lemma is a simple consequence of the maximum principle 
applied to equations (\ref{sR}) or (\ref{lR}). From
the lemma, we have the following relations
\[{\mathcal V}^s([g_0])\ge {\mathbb V}^s([g_0])\ge  {\mathcal V}^l([g_0])\ge {\mathbb V}^l([g_0]).\]
And we can show that the minimal volumes
${\cal V}^s([g_0])$ and ${\cal V}^l([g_0])$ are achieved.

\begin{cor}\label{vcomp} 
Suppose that
${\mathcal R}_{g_0}^s(x)< 0$ for any $x\in M$.
Then there is a unique conformal metric
$g^*\in [g_0]$ such that $vol({g^*})={\mathcal V}^s([g_0])$ with
${\mathcal R}_{g^*}^s(x)=-1, \forall x\in M$ and
\[{\mathbb V}^s([g_0])={\mathcal V}^s([g_0])\ge {\mathcal V}^l([g_0]).\]
The equality holds if and only if there is an Einstein metric in
$[g_0]$. If
${\mathcal R}_{g_0}^l(x)< 0$ for any $x\in M$, then there is a unique
$g^*\in [g_0]$ such that $vol(g_*)={\mathcal V}^l([g_0])$ with
${\mathcal R}_{g_*}^l(x)=-1, \forall x\in M$. In this case, we
have
\[{\mathcal V}^s([g_0])={\mathbb V}^s([g_0])\ge {\mathcal V}^l([g_0])={\mathbb V}^l([g_0]).\]
\end{cor}

For the study of minimal volumes in general Riemannian manifolds,
we refer to \cite{G}, \cite{BCG1} and \cite{Le}.

\section{positive curvature case}

Now we consider conformal classes with metrics of positive Ricci
curvature. Let $[g_0]$ be such a conformal class. After a suitable
scaling, we may assume that $Ric_{g_0}\ge (n-1)g_{0}$. Define
$[g_0]_{+}=\{g\in [g_0]\,|\, Ric_g\ge (n-1)g\}$ and
$V_{\max}(M,[g_0])=\sup_{g\in [g_0]_{+}} vol_g(M)$. This
definition is motivated by Gursky and Viaclovsky \cite{GV1}. From
the Bishop comparison, we know
\begin{equation}\label{eq5.0}
V_{\max}(M,[g_0]) \le vol(\S^n),\end{equation}
the volume of the unit sphere.

\begin{pro}\label{proY2} Let $(M,g_0)$ be a compact Riemannian manifold
with $Ric_{g_0}\ge (n-1)g_0$. If $V_{\max}(M,[g_0])<vol(\S^n)$,
then there is a conformal
metric
$g\in [g_0]_+$ with
\[{\mathcal R}^s_{g}=n-1.\] \end{pro}

\noindent{\it  Proof.} 
Consider the sequence of approximating function $F_k$ as in
Section 2 with a normalization condition that $F_k(1,1,\cdot,1)=n-1)$.
We first want to find a solution to the following
equation
\begin{equation}\label{eq5.1}
F_k\left(e^{2u}(\uuu)\right)=n-1,
\end{equation}
for large $k$. Define
\[V^k_{\max}([g_0])=\max\{ vol (g)\,|\, g \in [g_0] \hbox{ with }F_k(g^{-1}\cdot S_g
)\ge n-1\}.\] It is easy to check that
$\lim_{k\to \infty}V^k_{\max}([g_0])=V_{\max}([g_0]).$ Hence for large $k$ we have
\begin{equation}\label{eq5.3}
V^k_{\max}([g_0])<vol(\S^n).\end{equation}

To show the existence of solution of (\ref{eq5.1}), we consider a
deformation, which is similar to a deformation considered by
Gursky and Viaclovsky in their study of $\s_k$-Yamabe problem.
\begin{equation}\begin{array}{r}\label{eq5.2}
\ds\vs F_k\left(\n^2 u+du\otimes du-\frac 12 |\n u|^2 g+\psi(t)
S_g+(1-\psi)g\right)
\\ \ds =(n-1)(1-t) \left(\frac 1{vol(g_0)}\int_M e^{-(n+1)u}\right)^{\frac{2}{n+1}}
+(n-1)\psi(t) e^{-2u},\end{array}
\end{equation}
where $\psi(t): [0,1] \to [0,1]$ is a $C^1$ function satisfying
$\psi(0)=0$ and $\psi(t)=1$ for $t\ge 1/2$.
We now prove that there is a solution of (\ref{eq5.2})
when $t=1$, provided that $V^k_{\max}(M,[g])<vol(\S^n)$.

When $t=0$, it is easy to check that
(\ref{eq5.2}) has a unique solution
$u=0$ and its corresponding linearization has no nontrivial kernel.
Hence its Leray-Schauder degree  is non-zero.
If the solution space of (\ref{eq5.2}) for any
$t\in [0,1]$ is compact, then using degree theoretic argument,
we are done.
Assume by contradiction that there is no compactness.
Assume
without loss of generality that there is a sequence of solutions
$g_i=e^{-2u_i}g$ of (\ref{eq1}) with $t=1$ such that $u_i$ does
not converge in $C^{2,\a}$.
In view of Theorem \ref{thm1} and Corollary \ref{nthmC2}, we have
either
\begin{itemize}
\item[(a)] $\inf_{M}u_i \to -\infty$, or
\item[(b)] $\inf_{M}u_i \to +\infty$.
\end{itemize}
The latter is easy to be excluded as follows. At the minimum point
of $u_i$, we have
\[e^{2\inf u_i} F_k(S_{g_0}) =F_k(e^{2\inf u_i}S_0) \le (n-1),\]
which certainly implies that $\inf u_i$ is bounded from above uniformly.
Hence we are left to exclude (a).
  Let $x_i$ be the
minimum point of $u_i$ and assume that $x_i\to x_0$ as $i\to
\infty$.  Consider a scaled function
\[\tilde u_i=u(\exp_{x_i}\e_ix)-\log \e_i,\]
where $\e_i=\exp{u_i(x_i)}$. It is clear that $\tilde u_i\ge 0$
satisfies a similar equation on $B(0, \e_i^{-1} r_0/2) \subset
\R^n$ with a scaled metric, where $r_0$ is the injectivity radius
of $(M,g)$. The set $B(0, \e_i^{-1} r_0/2)$ with the scaled metric
converges to $\R^n$. By local estimates and local $C^2$ estimates
in \cite{GLW} for a more general concave case  one can show that
$\tilde u_i$ converges to an entire solution $v$ of
\begin{equation}
\label{eqa1.2} F_k(\n^2 v+dv\otimes dv-\frac 12|\n v|^2
g_{\R^n})=(n-1)e^{-2v}\end{equation} and
\begin{equation}
\label{eqa2.2}V^k_{\max}(g)\ge \lim_i\inf vol(e^{-2u_i}g)\ge vol(\R^n, e^{-2v} g_{\R^n}),
\end{equation}
where $g_{\R^n}$ is the standard Euclidean metric.
By a classification result of Li-Li in \cite{LL2},
we know that $(\R^n, e^{-2v} g_{\R^n})$ is equivalent to $\S^n$. Hence
$V^k_{\max}([g_0])\ge vol(\R^n, e^{-2v} g_{\R^n})= vol(\S^n)$.
 This  contradicts (\ref{eq5.3}).
Therefore, we have a solution $u_k$ of (\ref{eq5.2}) for large $k$.

Now we consider the sequence $\{u_k\}$. As above, we can show first that
$u_k$ has a uniform upper bound.
If $u_k$ has a uniform
lower bound,  Corollary \ref{nthmC2}
implies that the sequence $\{u_k\}$ has a uniform
$C^{2,\a}$ bound. And hence we have a limit $u_0$ which
is a solution we desire. Hence to prove the Proposition, we only
need to exclude the case that $\min u_k \to -\infty$. Assume that
we are in this case. By a similar argument presented above,
after considering a suitable rescaling we have  a limit $C^{2,\a}$
function $v_\infty$ satisfying
\begin{equation}
\label{eqa1} G_{n-1}(\n^2 v+dv\otimes dv-\frac 12|\n v|^2
g_{\R^n})=(n-1)e^{-2v}\end{equation} and
\begin{equation}
\label{eqa2}V^k_{\max}([g_0])\ge \lim_i\inf vol(e^{-2u_i}g)\ge
vol(\R^n, e^{-2v_\infty} g_{\R^n}).
\end{equation}
The contradiction follows from the
following Lemma. We finish the proof of the Proposition. \qed

\begin{pro}\label{proclass} Let $p<n/2$ and
$g=e^{-2u}g_{\R^n}$ be a $C^2$ function on $\R^n$ such that
\[G_p(g^{-1}S_g)(x)=c, \forall x\in \R^n,\]
for some constant $c$. Then $u=0$ if $c\le 0$ and $u(x)=\log
\frac{\l^2+|x-x_0|^2}{2\l \sqrt{\frac {(n-p)p}c }}$
if $c>0$. That is, $(\R^n, e^{-2u}g_{\R^n})$ can be compactified
as a standard sphere if $c>0$.
\end{pro}
\pr The proof follows \cite{LL2} closely. The only difference is
that the operator there is required to be $C^1$. Here our operator $G_p$
is Lipschitz only. However $G_p$ is uniformly elliptic and concave.
We will show in Lemma \ref{lem.a} below
that the Hopf lemma holds for our equation. Then the argument in \cite{LL2}
can be applied to our equation.

\qed

\begin{lem}\label{lem.a} Let $\O$ be a bounded domain in $\R^n$ and $p<n/2$.
If (\ref{eq3})
has  two solutions $w$ and $v$ with $w\ge v$ and $w(x_0)=v(x_0)$ for some
$x_0\in\partial \O$, then $x_0\in \partial \O$. Furthermore,
\[\frac {\partial w}{\partial \nu}(x_0) <\frac {\partial v}{\partial \nu}(x_0),\]
unless $w= v$. Here $\nu$ is the outer normal of $\partial \O$ at $x_0$.
\end{lem}
\pr For any function, set
\[A^u=\n u+du\otimes ds -\frac 12|\n u|^2g_0+S_{g_0}.\]
Since $G_p$ is concave and homogeneous one, we have
\[G_p(A^w-A^v)\le G_p(A^w)-G_p(A^v)=e^{-2w}-e^{-2v} \le 0.\]
Let $\tilde w=w-v$. Now we can  write $G_p(A^w-A^v)$ as follows
\[G_p(A^w-A^v)=a_{ij}(x)\tilde w_{ij}+b_i(x)\tilde w _i\]
with $\l_0 Id \le (a_{ij}(x)) \le \L_0 Id$ and $b_i(x)$ bounded
for any $i$. Therefore, we can apply Theorem 5 on page 61 and
Theorem 7 on page 65 in \cite{PW} to prove the Lemma.\qed

\begin{rem} Proposition \ref{proclass} does not hold for
$p=n-1$. For example $u=kx_1$ for any $k>0$ is a solution of
\[F(u)=G_{n-1}(\uuu)=0.\] The same example indicates that without the concavity  of $F$,
Proposition \ref{lem1} and Theorem \ref{thm1} are not true. It is
easy to check that $F(u_k)=0$. On a domain $\O\subset \{x_1\ge
0\}$, we have $u_k\ge 0$. But $|\n u_k|=k\to \infty$.  \end{rem}

\begin{pro}\label{proY1} Let $(M,g_0)$ be a compact Riemannian manifold
with $Ric_{g_0}\ge (n-1)g_0$. If $V_{\max}(M,[g_0])=vol(\S^n)$, then $(M, g)$
is conformally equivalent to the standard unit sphere. \end{pro}

This Proposition is a direct consequence of the following
\begin{pro}\label{pro5.2} Let $(M,g)$ be a compact Riemannian manifold with
$Ric_g\ge n-1$. If  $vol(M)$ is close to $\o_n$, the volume of
$\S^n$, then, the Yamabe constant of $(M,[g])$, $Y(M,[g])$ is close
to  $n(n-1) \o_n^{2/n}$, the  Yamabe constant of $\S^n$.
\end{pro}

 \pr Let us first recall the well-known Yamabe constant
of $(M,g)$, which  is defined by
\[Y(M,[g]):=\inf \left(\int  v^{\frac{2n}{n-2}} dvol(g)\right)^{-\frac{n-2}n}
 \left\{4\frac{n-1}{n-2}\int|\n v|^2 dvol(g)+\int R_g v^2 dvol(g)\right\}
 \]
By a result of Ilias \cite{I}, which is based on a result of Gromov
(see also \cite{BCG0}), we have
 \[ n(n-1)\o_n^{2/n}\left(\int v^{\frac{2n}{n-2}}\right)^{(n-2)/2}\le
 \left(\frac{\o_n}{vol(g)}\right)^{2/n}
 \left\{4\frac{n-1}{n-2}\int|\n v|^2 +n(n-1)\int  v^2
 \right\},\]
 for any $v\in H^{2}_1(M)$. Note that $R_g\ge n(n-1)$. Therefore, we have
\[Y(M,[g])\ge  \left(\frac{\o_n}{vol(g)}\right)^{-2/n}n(n-1)\o_n^{2/n}
\ge n(n-1)\o_n^{2/n}-\d,\] for any small $\d>0$, provided that $
vol(g)$ is close to $\o_n$. \qed

By the results of Colding (see \cite{colding1},\cite{colding2} and
\cite{petersen}), we know that the condition in Proposition
\ref{pro5.2} is equivalent to the one of following three other conditions:
\begin{itemize}
\item [1)] rad$M$ is close to $\pi$,
\item [2)] $M$ is Gromov-Hausdorff close to $\S^n$,
\item [3)] the $(n+1)$th eigenvalue  of the Laplacian,
 $\l_{n+1}(M)$, is close to $n$.
\end{itemize}

\

\noindent{\it Proof of Proposition \ref{proY1}.}
$V_{\max}(M,[g_0])=vol(\S^n)$ implies by definition that there is
a sequence $g_i\in [g_0]_+$ with $\lim_{i\to \infty}
vol(g_i)=vol(\S^n)$. Proposition \ref{pro5.2} implies that
$Y(M,[g_0])=Y(M,[g_i]) \to n(n-1)\o_{n}^{2/n}$, the Yamabe
constant of $\S^n$. Hence, the  Yamabe constant of $(M,[g_0])$
equals to  the Yamabe constant of the standard sphere. By the
resolution of the Yamabe problem by Aubin \cite{Aubin} and Schoen
\cite{schoen1}, $(M,[g_0])$
is conformally equivalent to the standard sphere. \qed

\

\noindent{\it Proof of Theorem \ref{thm3}.} Theorem \ref{thm3}
follows from Propositions \ref{proY1} and \ref{proY2}. \qed

\begin{rem} It is interesting to know weather $V_{\max}$ is achieved as
in the negative case. One can show that if $V_{\max}$ is achieved
by $\tilde g$, then ${\mathcal R}^s_{\tilde g}$ is constant.
\end{rem}

Now we prove Theorem \ref{thm5}.

\noindent{\it Proof of Theorem \ref{thm5}.} It follows the exact
same arguments in the proof of Theorem 3 in \cite{GLW}, since that
proof works for general uniform elliptic concave equations as
well, as we note that $W_p(g)>0$ implies the positivity of the
mean curvature when $p\le \frac n2$. We only give a sketch here.

\

\noindent{\bf Step 1.} We define a deformation
\begin{equation}\label{defo}
f_t(g):=tW_p(g)+(1-t)R_g=1,\end{equation}
where $g=e^{-2u}g_0$. Equation (\ref{defo}) ($\forall t \in [0,1]$)
is still uniformly elliptic and concave.

\

\noindent{\bf Step 2.}  (Harnack inequality) There is a constant
$C>0$ such that for  a solution $u$ of (\ref{defo})
in $B_{3R}$ we have
\begin{equation}\label{har}
\min_{B_R}u+\max_{B_{2R}}u\ge 2\log R-\log C.\end{equation}
Here $B_R$ is the ball of radius $R$ in $\R^n$.
(\ref{har}) can be proved as in \cite{GLW} using the method of
moving planes.

By scaling argument, we may assume that $R=1$. Assume by
contradiction that (\ref{har}) is not true. Then there exists a
sequence of solutions of (\ref{defo}) in $B_3$ such that
 \begin{equation}\label{x1}
 \min_{B_1} u_i+ \max_{B_2} u_i <- i.\end{equation}
Let $m_i=  \min_{B_1} u_i=u_i(\bar x_i)$, and let $x_i\in B_1$
with $\bar B_{r_i}(x_i)\subset \bar B_1$ and $|x_i-\bar x_i|=r_i$.
Here $r_i=e^{m_i}$. In view of (\ref{x1}), we know that $r_i\to 0$
as $i \to \infty$. Consider  a new sequence of functions $v_i$
defined by
 \[\begin{array}{rcl}
 v_i(x)&=& u_i(x_i+r_ix)-m_i\end{array}\]
and set $\bar x_i=x_i+r_i\bar y_i.$ It is clear that $v_i$
satisfies (\ref{defo}) in $\{|x|<r_i^{-1}\}$ and $|\bar y_i|=1$.
>From (\ref{x1}), we extend $v_i$ by the Kelvin transformation to
the whole Euclidean space as in (2.8) of \cite{GLW}. Now applying
the method of moving planes as in \cite{GLW}, which in turn
follows closely from \cite{chenLin}, we can show that $v_i$
converges to $0$ in $B_{\frac 12}(\bar y_i)$. This is a
contradiction. Note that though we are dealing with the Lipschitz
operators, the method of moving planes works by using the fact
that $f_t$ in (\ref{defo}) is uniformly elliptic and concave.

\

\noindent{\bf Step 3.} Consider a solution $u$ of equation
(\ref{defo}). First, it is clear that we have that the scalar
curvature of $g=e^{-2u}g_0$ is positive. Hence we can apply the
result of Schoen-Yau in \cite{SchYau} to embed the universal cover
$\widetilde M$ of $(M,g)$ into $\S^n$ by a map $\Phi$ conformally.
Therefore we can use the method of moving planes (again make use
of uniformly ellipticity and concavity of $f_t$) to obtain as in
\cite{GLW} that
 \[|\n u|(x)\le C, \quad \hbox{ for any } x\in M,\]
for some constant independent of $u$, provided that $(M,g_0)$ is
not equivalent to $\S^n$. It follows that
\begin{equation}\label{est}
\max u-\min u \le C,\end{equation} for some constant independent
of $u$. (\ref{est}), together with the Harnack
inequality (\ref{har}), implies that
\[\min u \ge C,\]
for some constant independent of $u$. Hence by Theorem
\ref{thmC2}, we know that the solution space of equation
(\ref{defo}) is compact.

\

\noindent{\bf Step 4.} From Step 3 we can apply the degree theory.
We may use a result of Li in \cite{Li89}, a variation of the
original Leray-Schauder theorem \cite{LS}. We also refer to
Nirenberg's lecture notes \cite{Nirenberg} on the exposition of
the degree theory in nonlinear differential equations. When $t=0$,
the topological degree for equation (\ref{defo}) is $-1$, which
was proved by Schoen \cite{18}. Since the solution space is
compact, the topological degree for equation (\ref{defo}) with
$t=1$ is also $-1$. This finishes the proof of the Theorem.

\qed

\begin{rem} As in Remark \ref{rm1}, the conditions in Theorem \ref{thm3}
and Theorem \ref{thm5} can be weakened to the assumption that the
corresponding curvature of the background metric is  nonnegative
and positive at some point. The same argument using the
short time existence of the corresponding curvature flows as in
Remark \ref{rm1} can produce a metric $g\in [g_0]$ with the
positive curvatures.
\end{rem}

\section{Appendix}
 In this Appendix, we assume that $(M,g)$ is a compact \lcf manifold.
 Recall the
Weitzenb\" ock formula for $p$-forms $\omega$
\[\D \omega  =\n^* \n  \omega +{\cal R}\omega,\]
where
\[{\cal R}\omega =\sum_{j,l=1} \omega _j \wedge i(e_l) R(e_j,e_l) \omega.\]
Here $e_j$ is a local basis and $i(\cdot)$ denotes the interior
product $\D=dd^*+d^*d$ is the Hodge-de Rham Laplacian and $\n^*
\n$ is the (positive) Laplacian.
 In local coordinates, let
$\omega =\omega_{1} \wedge \cdots \wedge \omega_{p}$.  Since $(M,g)$ is
locally conformally flat, we have 
\begin{equation}\label{Weiz}{\cal R} \omega
=\left((n-p) \sum_{i=1}^p \l_i+ p\sum_{i=p+1}^n
\l_i\right)\omega,\end{equation} where $\l$'s are eigenvalues of
the Schouten tenser $S_g$. See for instance \cite{N} or \cite{GLW}.
By the Bochner technique, we have 
\begin{pro} Let $(M,g)$ is a compact \lcf manifold and $1\le p\le n/2$.
If  its $p$-Weitzenb\"ock curvature $W_p(g)$ is positive,
then the Betti number $b_q=0$ for $p\le q \le n-q$.
If its $p$-Weitzenb\"ock curvature $W_p(g)$ is non-negative and $b_p\neq 0$,
then $M$ is a quotient of $H^p\times \S^{n-p}$ if $p>1$
and a quotient of $\S^1\times \S^{n-1}$ if $p=1$, where $H^p$ is the hyperbolic
space of curvature $-1$ and $\S^{n-p}$ is the standard sphere of  curvature $1$. 
\end{pro}

For the vanishing of cohomology groups under various conditions,
see \cite{N} and \cite{GLW} and referrences therein.

\end{document}